\documentclass[12pt]{article}
\usepackage[top=1in,bottom=1in,left=1in,right=1in]{geometry}
\usepackage{indentfirst}
\usepackage{amsfonts,amsmath,amsthm,amssymb,hyperref,enumerate,bm}
\usepackage{longtable}
\usepackage[caption=false]{subfig}
\usepackage{leftidx}
\usepackage{lineno}
\usepackage{color,xcolor}
\usepackage{array}
\usepackage{latexsym}
\usepackage{float}
\usepackage{color}
\usepackage{longtable,supertabular,booktabs}
\usepackage{rotating}
\usepackage{cases}
\usepackage{multicol,multirow}
\usepackage{epsfig}
\usepackage{cite}
\usepackage{fancyhdr}       % fuer Deckblatt
\usepackage{dsfont}
\usepackage{hyperref}
\hypersetup{
	colorlinks=true,
	linkcolor=red,
	filecolor=blue,
	urlcolor=red,
	citecolor=blue,
}

\newtheorem{theorem}{Theorem}[section]

\newtheorem{remark}{Remark}[section]

\def\G1{G^\mathcal{C}}

 \newenvironment{prof}{\trivlist
      \item[\hskip\labelsep
      {\itshape Proof.}]\normalfont}
      {\hspace*{\fill}$\Box$\endtrivlist}

\begin{document}
\title{$(G,F)$-points on $\mathbb{Q}$-algebraic varieties}
\author{
Yangcheng Li$^{1,}$\footnote{E-mail\,$:$ liyangchengmlx@163.com.}\quad\,
Hongjian Li$^{2,}$\footnote{Corresponding author. E-mail\,$:$ lhj@gdufs.edu.cn. Supported by the National Natural Science Foundation of China (Grant No. 12171163), the Basic and Applied Basic Research Foundation of Guangdong Province (Grant No. 2024A1515010589) and the Project of Guangdong University of Foreign Studies (Grant No. 2024RC063).}\\
{\small\it  $^{1}$School of Mathematical Sciences, South China Normal University,}\\
{\small\it Guangzhou 510631, Guangdong, P. R. China} \\
{\small\it  $^{2}$School of Mathematics and Statistics, Guangdong University of Foreign Studies,}\\
{\small\it Guangzhou 510006, Guangdong, P. R. China} \\
}

\date{}
 \maketitle
\date{}

\noindent{\bf Abstract}\quad
Let \(G\in \mathbb{Q}[x,y,z]\) be a polynomial, and let $V(G)$ be the $\mathbb{Q}$-algebraic variety corresponding to $G$, i.e., $V(G)=\{P\in\mathbb{Q}^3~|~G(P)=0\}$. Let
\[\begin{split}
F:\quad  &\mathbb{Q}^3\rightarrow \mathbb{Q}^3,\\
&(x,y,z)\mapsto (f(x),f(y),f(z))
\end{split}\]
be a vector function, where $f\in \mathbb{Q}[x]$. It is easy to know that the function obtained by the composition of \(G\) and \(F\), denoted as \(G\circ F\), is still in \(\mathbb{Q}[x,y,z]\). Moreover, let $V(G\circ F)$ be the $\mathbb{Q}$-algebraic variety corresponding to $G\circ F$, i.e., $V(G\circ F)=\{P\in\mathbb{Q}^3~|~G\circ F(P)=0\}$. A rational point \(P\) is called a \((G,F)\)-point on \(V(G)\) if \(P\) belongs to the intersection of \(V(G)\) and \(V(G\circ F)\), that is \(P\in V(G)\cap V(G\circ F)\). Denote \(\langle G,F\rangle\) as the set consisting of all \((G,F)\)-points on \(V(G)\). Obviously, \(\langle G,F\rangle\) is a \(\mathbb{Q}\)-algebraic variety. In this paper, we consider the algebraic variety \(\langle G,F\rangle\) for some specific functions \(G\) and \(F\). For these specific functions $G$ and $F$, we prove that \(\langle G,F\rangle\) will be isomorphic to a certain elliptic curve. We also analyze some properties of these elliptic curves.


\medskip
\noindent{\bf MR(2020) Subject Classification} 11D25, 11D72, 14A10, 11G05.

\section{Introduction}
Diophantine geometry aims to study the rational points on algebraic varieties. However, it is not an easy task to determine the rational points on an algebraic variety. In other words, it is not easy to find the rational solutions of certain Diophantine equations. The famous Mordell Conjecture states that there are only finitely many rational points on a curve with genus greater than $1$. In 1983, the Mordell Conjecture was finally proven by Faltings. Curves with genus equal to $0$ can be parameterized. Therefore, we are concerned with curves with genus equal to $1$. The representative of curves with genus equal to $1$ is the elliptic curve. However, the structure of rational points on elliptic curves is still not fully understood.

Several authors (see \cite{Bennett,Katayama,Schinzel-Sierpinski,Ulas1,Ulas2,Zhang-Cai1,Zhang-Cai2,Zhang1,Zhang2,Zhang-Zargar}) investigated the rational solutions of the Diophantine equation
\begin{equation*}
f(x)f(y)=f(z)^2,
\end{equation*}
where $f\in \mathbb{Q}[x]$. In this studies, the transformation
\[x=T,~~y=u^2T,~~z=uT\]
was frequently used. It means that they actually gave the rational solutions of the Diophantine system
\begin{equation}
\begin{cases}
\begin{aligned}
&f(x)f(y)=f(z)^2,\\
&xy=z^2.                                                            \label{1.1}
\end{aligned}
\end{cases}
\end{equation}
Naturally, we can consider the rational solutions of the general Diophantine system
\begin{equation}
\begin{cases}
\begin{aligned}
&G(f(x),f(y),f(z))=0,\\
&G(x,y,z)=0,                                                            \label{1.2}
\end{aligned}
\end{cases}
\end{equation}
where \(G\in \mathbb{Q}[x,y,z]\) and $f\in \mathbb{Q}[x]$. Let us redescribe the Diophantine system (\ref{1.2}). Let \(G\in \mathbb{Q}[x,y,z]\) be a polynomial, and let $V(G)$ be the $\mathbb{Q}$-algebraic variety corresponding to $G$, i.e., $V(G)=\{P\in\mathbb{Q}^3~|~G(P)=0\}$. Let
\[\begin{split}
	F:\quad  &\mathbb{Q}^3\rightarrow \mathbb{Q}^3,\\
	&(x,y,z)\mapsto (f(x),f(y),f(z))
\end{split}\]
be a vector function, where $f\in \mathbb{Q}[x]$. It is easy to know that the function obtained by the composition of \(G\) and \(F\), denoted as \(G\circ F\), is still in \(\mathbb{Q}[x,y,z]\). Moreover, let $V(G\circ F)$ be the $\mathbb{Q}$-algebraic variety corresponding to $G\circ F$, i.e., $V(G\circ F)=\{P\in\mathbb{Q}^3~|~G\circ F(P)=0\}$. A rational point \(P\) is called a \((G,F)\)-point on \(V(G)\) if \(P\) belongs to the intersection of \(V(G)\) and \(V(G\circ F)\), that is \(P\in V(G)\cap V(G\circ F)\). Denote \(\langle G,F\rangle\) as the set consisting of all \((G,F)\)-points on \(V(G)\). Obviously, \(\langle G,F\rangle\) is a \(\mathbb{Q}\)-algebraic variety.

When $G(x,y,z)=xy-z^2$, \(\langle G,F\rangle\) represents all the rational solutions of the Diophantine system (\ref{1.1}). From \(G(x,y,z)=xy-z^2=0\), we obtain \(z=\sqrt{xy}\), which means that \(z\) is the geometric mean of \(x\) and \(y\).

For some specific functions \(G\) and \(F\), the \(Q\)-algebraic variety $V(G\circ F)$ may be an empty set or may have only trivial solutions. When \(G = x + y - z\) and \(f(x)=x^{n}\), by Fermat's Last Theorem, $V(G\circ F)$ has only trivial solutions. When \(G = xy - z^{2}\) and \(f(x)=x^{k}-1\), Bennett \cite{Bennett} proved that $V(G\circ F)$ has only trivial solutions. When \(G = x^{2}+y^{2}-z^{2}\) and \(f(x)=x(x + 1)/2\), Sierpi\'nski \cite{Sierpinski} introduced (given by Zarankiewicz) a nontrivial positive integer solution \((x, y,z)=(132,143,164)\) for $V(G\circ F)$.

In 2019, Zhang and Chen \cite{Zhang-Chen} studied the rational solutions of Diophantine equation of harmonic mean
\[f(z)=\frac{2}{\frac{1}{f(x)}+\frac{1}{f(y)}}.\]

In this paper, we consider the algebraic variety \(\langle G,F\rangle\) for some specific functions \(G\) and \(F\). For the function $G$, we consider that \(G=xy-z^{2}\), or \(G=(x+y)z-2xy\). For \(G=(x+y)z-2xy=0\), it means that $z$ is the harmonic mean of $x$ and $y$. For the vector function $F$, we consider the components $f$ of $F$ to be \(f=ax^{2}+bx+c\), \(f=ax+b+cx^{-1}\), \(f=x(ax^{2}+bx+c)\), respectively. For these specific functions $G$ and $F$, we prove that \(\langle G,F\rangle\) will be isomorphic to a certain elliptic curve. We also analyze some properties of these elliptic curves.

\section{The main results}
By the theory of elliptic curves, we prove
\begin{theorem}
	Let \(G = xy - z^{2}\) and \(f(x) = ax^{2}+bx + c \in \mathbb{Q}[x]\) with \(abc \neq 0\).
When \(4ac - b^{2}\neq 0\), \(\langle G,F\rangle\) is birationally equivalent to the elliptic curve
\[\mathcal{E}_1: Y^{2}=X^3+27a^2b^2c^2(3ac-b^2)X+27a^3b^4c^3(9ac-2b^2).\]
When \(4ac - b^{2}=0\), \(\langle G,F\rangle\) is a curve with genus $0$ and its parameterization is given by
\[y=-\frac{bct^2}{(2ct + b)^2},~~z=-\frac{(bt + 4ct + 2b)ct}{(2ct + b)(bt + 2ct + b)},\]
where $t$ is a rational number.
\end{theorem}

\begin{prof}
When \(G=xy-z^{2}\) and \(f(x)=ax^{2}+bx+c\), the Diophantine system (\ref{1.2}) is equivalent to
\begin{equation}
	(y-z)^2\left(abyz^2+acy^2+2acyz+acz^2+bcy\right)=0.               \label{2.1}
\end{equation}
Since \(y=z\) is trivial, we only need to consider the rational points on the curve
\[\mathcal{C}_1: abyz^2+acy^2+2acyz+acz^2+bcy=0.\]
By the map $\varphi_1$:
\begin{equation}
	\begin{split}
	X=&~\frac{3abc(3acy+6acz-b^2y+2bc)}{by+c},\\
	Y=&~\frac{27a^2bc^2(abcy^2+3abcyz-b^3yz-ac^2y-ac^2z+b^2cy-b^2cz-bc^2)}{(by+c)^2},       \label{2.2}
	\end{split}
\end{equation}
with the inverse map $\varphi^{-1}_1$:
\begin{equation}
	\begin{split}
	y=&~-\frac{c(108a^3b^2c^3-18a^2b^4c^2+9Xa^2c^2-3Xab^2c+6Yac+X^2)}{b(X-9a^2c^2+3ab^2c)^2},\\
	z=&~-\frac{9a^2b^2c^2+3Xac+Y}{3ba(X-9a^2c^2+3ab^2c)},                                                \label{2.3}
	\end{split}
\end{equation}
we can transform $\mathcal{C}_1$ into the elliptic curve
\begin{align*}
	\mathcal{E}_1:~Y^{2}=X^3+27a^2b^2c^2(3ac-b^2)X+27a^3b^4c^3(9ac-2b^2).
\end{align*}
The discriminant of $\mathcal{E}_1$ is
\[\Delta_1(E)=-531441a^8b^6c^8(4ac-b^2).\]
Hence, if $4ac-b^2\neq0$, then $\Delta_1(E)\neq0$, so $\mathcal{E}_1$ is non-singular. Therefore, (\ref{2.2}) and (\ref{2.3}) give a bijection between \(\langle G,F\rangle\) and the elliptic curve $\mathcal{E}_1$, so \(\langle G,F\rangle\) is birationally equivalent to the elliptic curve $\mathcal{E}_1$.

When $4ac-b^2=0$, the Diophantine system (\ref{1.2}) is equivalent to
\begin{equation}
	b^2yz^2 + bcy^2 + 2bcyz + bcz^2 + 4c^2y=0.                                    \label{2.4}
\end{equation}
The curve given by (\ref{2.4}) is a curve of genus $0$, and its parameterization is given by
\[y=-\frac{bct^2}{(2ct + b)^2},~~z=-\frac{(bt + 4ct + 2b)ct}{(2ct + b)(bt + 2ct + b)},\]
where $t$ is a rational number. This completes the proof.
\end{prof}

\begin{remark}
Since \(xy=z^{2}\), when \(f=x^k(ax^2+bx+c),k\in\mathbb{Z}\), we can obtain exactly the same result as that of Theorem 1.1. Particularly, when \(k = 1\) or \(k = -1\), we obtain \(f=x(ax^2+bx+c)\) and \(f=ax+b+cx^{-1}\), respectively.
\end{remark}

\begin{theorem}
When \(b^{2}\neq kac~(k=-1,1,3)\), the elliptic curve \(\mathcal{E}_1\) has a positive rank. When \(b^{2}=kac~(k=-1,1,3)\), the rank of the elliptic curve \(\mathcal{E}_1\) is \(0\).
\end{theorem}

\begin{prof}
It is easy to check that the elliptic curve $\mathcal{E}_1$ contains two rational points
\begin{align*}
	P_0=(-3ab^2c,0),~~P_1=(6ab^2c, 27a^2b^2c^2).
\end{align*}
By the group law, we have
\begin{align*}
	P_2=[2]P_1=\bigg(&\frac{3(3ac-b^2)(ac-3b^2)}{4},~-\frac{27(ac-b^2)(a^2c^2+4ab^2c-b^4)}{8}\bigg),\\
	P_3=[3]P_1=\bigg(&\frac{6ab^2c(13a^4c^4+24a^3b^2c^3-22a^2b^4c^2+b^8)}{(a^2c^2-6ab^2c+b^4)^2},\\
	                 &-\frac{27a^2b^2c^2(3ac-b^2)(ac+b^2)(a^4c^4+24a^3b^2c^3-22a^2b^4c^2+16ab^6c-3b^8)}{(a^2c^2-6ab^2c+b^4)^3}\bigg)\\
\end{align*}
and
\begin{align*}
	P_4&=[4]P_1\\
       &=\bigg(\frac{3(3ac-b^2)X_4}{16(ac-b^2)^2(a^2c^2+4ab^2c-b^4)^2},\frac{27(a^4c^4-20a^3b^2c^3+6a^2b^4c^2-4ab^6c+b^8)Y_4}{64(ac-b^2)^3(a^2c^2+4ab^2c-b^4)^3}\bigg),
\end{align*}
where
\begin{align*}
	X_4=&~a^7c^7-45a^6b^2c^6+365a^5b^4c^5-121a^4b^6c^4\\
	&+307a^3b^8c^3-151a^2b^{10}c^2+31ab^{12}c-3b^{14},\\
	Y_4=&~a^8c^8+80a^7b^2c^7-180a^6b^4c^6+656a^5b^6c^5\\
	&-282a^4b^8c^4-80a^3b^{10}c^3+ 76a^2b^{12}c^2-16ab^{14}c+b^{16}.
\end{align*}
Let the line go through the points $P_0$ and $P_1$, intersecting $\mathcal{E}_1$ at $P_5$, then
\begin{align*}
	P_5&=-(P_0+P_1)=\left(3ac(3ac-b^2),27a^3c^3~\right).
\end{align*} Similarly,
\begin{align*}
	P_6=&-(P_0+P_2)=\left(\frac{3ab^2c(11a^2c^2+2ab^2c-b^4)}{(ac-b^2)^2},-\frac{54(a^2c^2+4ab^2c -b^4)a^3b^2c^3}{(ac-b^2)^3}\right),\\
	P_7=&-(P_0+P_3)=\bigg(\frac{3ac(3a^5c^5-45a^4b^2c^4+102a^3b^4c^3-34a^2b^6c^2+7ab^8c-b^{10})}{(ac+b^2)^2(3ac-b^2)^2},\\
	&-\frac{27c^3a^3(a^4c^4+24a^3b^2c^3-22a^2b^4c^2+16ab^6c-3b^8)(a^2c^2-6ab^2c+b^4)}{(3ac-b^2)^3(ac+b^2)^3}\bigg),\\
\end{align*}
and
\begin{align*}
	P_8=-(P_0+P_4)=\bigg(&\frac{3ab^2cX_8}{(a^4c^4-20a^3b^2c^3+6a^2b^4c^2-4ab^6c+b^8)^2},\\
	&~\frac{108c^3b^2a^3(ac-b^2)(a^2c^2+4ab^2c-b^4)Y_8}{(a^4c^4-20a^3b^2c^3+6a^2b^4c^2-4ab^6c+b^8)^3}\bigg),
\end{align*}
where
\begin{align*}
	X_8=&~47a^8c^8+328a^7b^2c^7-460a^6b^4c^6-1096a^5b^6c^5\\
	&+1290a^4b^8c^4-392a^3b^{10}c^3+20a^2b^{12}c^2+8ab^{14}c- b^{16},\\
	Y_8=&~a^8c^8+80a^7b^2c^7-180a^6b^4c^6+656a^5b^6c^5\\
	&-282a^4b^8c^4-80a^3b^{10}c^3+76a^2b^{12}c^2 -16ab^{14}c+b^{16}.
\end{align*}

Let $X(P)$ denote the $X$-coordinate of the point $P$. When $a\neq 0$, we give the conditions such that $X(P_i)=X(P_j),~0\leq i<j\leq8$ in the following table.
\begin{table}[h]
\centering
	\begin{tabular}{ll|l}
		\toprule
		&$X$-coordinate & $k~(b^2=kab)$   \\
		\midrule
		&$X(P_0)=X(P_i),~i=1,...,8$ &$-1,~1,~3$  \\
		&$X(P_1)=X(P_i),~i=2,...,8$ &$-1,~1,~3$  \\
		&$X(P_2)=X(P_i),~i=3,...,8$ &$-1,~1,~3$  \\
		&$X(P_3)=X(P_i),~i=4,...,8$ &$1$   \\
		&$X(P_4)=X(P_i),~i=5,...,8$  &$-1,~3$   \\
		&$X(P_5)=X(P_i),~i=6,...,8$  &$1$  \\
		&$X(P_6)=X(P_i),~i=7,8$  &$-1,~3$ \\
		&$X(P_7)=X(P_8),~i=8$         & none   \\
		\bottomrule
	\end{tabular}
	\vspace{0.23cm}  \caption{\label{tab:6} Conditions such that $X(P_i)=X(P_j), 0\leq i<j\leq8$.}
\end{table}

Table 1 shows that when \(b^{2}\neq kac~(k=-1,1,3)\), the points $P_0$ and $\pm P_i,~i=1,...,8,$ are different. By Mazur's theorem (see p. 58 of \cite{Silverman-Tate}) about the quantity of rational points and the rank of elliptic curve: If an elliptic curve $E$ defined over $\mathbb{Q}$ has more than 16 different rational points, then it has infinitely many rational points and its rank has at least one. Therefore, $\mathcal{E}_1$ has a positive rank, and thus there are infinitely many rational points on $\mathcal{E}_1$.

When \(b^{2}=kac~(k=-1,1,3)\), we have \(a=\frac{b^2}{ck}\), and then
\begin{align*}
	\mathcal{E}'_1:~Y^2=~X^3-\frac{27(k-3)b^8}{k^3}X-\frac{27(2k-9)b^{12}}{k^4}.
\end{align*}
By the transformation
\begin{equation}
		U=\frac{Y}{b^6},\quad V=\frac{X}{b^4},
\end{equation}
we get
\begin{align*}
	\mathcal{E}_{(k)}:~U^2=V^3-\frac{27(k-3)}{k^3}V-\frac{27(2k-9)}{k^4}.
\end{align*}

1) When $k=-1$, we get
\begin{align*}
	\mathcal{E}_{(-1)}:~U^2=V^3-108V+297.
\end{align*}
Using the package of Magma, the rank of $\mathcal{E}_{(-1)}$ is $0$, and the only rational points on $\mathcal{E}_{(-1)}$ are
\[(3,0),~(-6, \pm 27)~(12,\pm 27).\]

2) When $k=1$, we get
\begin{align*}
	\mathcal{E}_{(1)}:~U^2=V^3+54V+189.
\end{align*}
The rank of $\mathcal{E}_{(1)}$ is $0$, and the only rational points on $\mathcal{E}_{(1)}$ are
\[(-3,0),~(6, \pm 27).\]

3) When $k=3$, we get
\begin{align*}
	\mathcal{E}_{(3)}:~U^2=V^3+1.
\end{align*}
The rank of $\mathcal{E}_{(3)}$ is $0$, and the only rational points on $\mathcal{E}_{(3)}$ are
\[(-1,0),~(0,\pm 1),~(2,\pm3).\]

Therefore, when \(b^{2}\neq kac~(k=-1,1,3)\), the elliptic curve \(\mathcal{E}_1\) has a positive rank. When \(b^{2}=kac~(k=-1,1,3)\), the rank of the elliptic curve \(\mathcal{E}_1\) is \(0\).
\end{prof}

\begin{theorem}
	Let \(G=(x+y)z-2xy\) and \(f(x)=ax^{2}+bx+c \in \mathbb{Q}[x]\) with \(abc \neq 0\). When \(2ac - b^{2}\neq 0\) and \(4ac-b^2\neq0\), \(\langle G,F\rangle\) is birationally equivalent to the elliptic curve
	\[\mathcal{E}_2: Y^{2}=X^3-3a^4c^4X+a^4c^4(2a^2c^2-4ab^2c+b^4).\]
When $2ac - b^{2}= 0$, \(\langle G,F\rangle\) is the union of two curves, which are
\[z=-\frac{2c}{b},~~z=\frac{2cy(by + 2c)}{b^2y^2 + 2c^2}.\]
When $4ac - b^{2}= 0$, \(\langle G,F\rangle\) is a curve with genus $0$ and its parameterization is given by
\[\begin{split}
y=&~\frac{24ct(2bt - 12ct - c)}{68b^2t^2 + 48bct^2 - 144c^2t^2 + 4bct - 24c^2t - c^2},\\
z=&-\frac{48ct(8bt - 12ct - c)}{100b^2t^2 - 192bct^2 + 144c^2t^2 - 16bct + 24c^2t + c^2},
\end{split}\]
where $t$ is a rational number.
\end{theorem}

\begin{prof}
When \(G=(x+y)z-2xy\) and \(f(x)=ax^{2}+bx+c\), the Diophantine system (\ref{1.2}) is equivalent to
\begin{equation}
	(y-z)^2\left(a^2y^2z^2-2acy^2-2acyz+acz^2-2bcy+bcz\right)=0.               \label{3.1}
\end{equation}
Since \(y=z\) is trivial, we only need to consider the rational points on the curve
\[\mathcal{C}_2: a^2y^2z^2-2acy^2-2acyz+acz^2-2bcy+bcz=0.\]
By the map $\varphi_2$:
\begin{equation}
	\begin{split}
		X=&~\frac{c(a^2by^2z+a^2cy^2+abcz+b^2c)}{y^2},\\
		Y=&~\frac{bc^2(a^3y^3z+a^2by^3+a^2by^2z+2a^2cy^2+a^2cyz+abcy+abcz+b^2c)}{y^3},       \label{3.2}
	\end{split}
\end{equation}
with the inverse map $\varphi^{-1}_2$:
\begin{equation}
	\begin{split}
		y=\frac{bc(-a^3c^3+a^2b^2c^2+Xac+Y)}{a^4c^4-2a^3b^2c^3-2Xa^2c^2+X^2},~~z=\frac{2bc(-2a^3c^3-a^2b^2c^2+2Xac+Y)}{a^4c^4+4a^3b^2c^3-2Xa^2c^2+X^2},                                                \label{3.3}
	\end{split}
\end{equation}
we can transform $\mathcal{C}_2$ into the elliptic curve
\begin{align*}
	\mathcal{E}_2:~Y^{2}=X^3-3a^4c^4X+a^4c^4(2a^2c^2-4ab^2c+b^4).
\end{align*}
The discriminant of $\mathcal{E}_2$ is
\[\Delta_2(E)=27a^8c^8b^2(4ac-b^2)(2ac-b^2)^2.\]
Hence, if $2ac-b^2\neq0$ and $4ac-b^2\neq0$, then $\Delta_2(E)\neq0$, so $\mathcal{E}_2$ is non-singular. Therefore, (\ref{3.2}) and (\ref{3.3}) give a bijection between \(\langle G,F\rangle\) and the elliptic curve $\mathcal{E}_2$, so \(\langle G,F\rangle\) is birationally equivalent to the elliptic curve $\mathcal{E}_2$.

When $2ac-b^2=0$, the Diophantine system (\ref{1.2}) is equivalent to
\begin{equation*}
	(bz + 2c)(b^2y^2z - 2bcy^2 - 4c^2y + 2c^2z)=0,
\end{equation*}
which leads to
\[z=-\frac{2c}{b},~~z=\frac{2cy(by + 2c)}{b^2y^2 + 2c^2}.\]
When $4ac-b^2=0$, the Diophantine system (\ref{1.2}) is equivalent to
\begin{equation}
	b^3y^2z^2 - 8bc^2y^2 - 8bc^2yz + 4bc^2z^2 - 32c^3y + 16c^3z=0.                                    \label{3.4}
\end{equation}
The curve given by (\ref{3.4}) is a curve of genus $0$, and its parameterization is given by
\[\begin{split}
y=&~\frac{24ct(2bt - 12ct - c)}{68b^2t^2 + 48bct^2 - 144c^2t^2 + 4bct - 24c^2t - c^2},\\
z=&-\frac{48ct(8bt - 12ct - c)}{100b^2t^2 - 192bct^2 + 144c^2t^2 - 16bct + 24c^2t + c^2},
\end{split}\]
where $t$ is a rational number. This completes the proof.
\end{prof}

\begin{theorem}
The elliptic curve \(\mathcal{E}_2\) has a positive rank.
\end{theorem}

\begin{prof}
It is easy to check that the elliptic curve $\mathcal{E}_2$ contains the following rational point
\begin{align*}
	P=\left(-a^2c^2, -c^2(2ac-b^2)a^2\right).
\end{align*}
By the group law, we get the following eight points
	\begin{align*}
		&[2]P,\quad [3]P,\quad [4]P,\quad [5]P,\quad [6]P,\quad [7]P,\quad [8]P,\quad [9]P.
	\end{align*}
	The points $[2]P$ and $[3]P$ are as follows
\begin{equation*}
	\begin{split}
	[2]P=&\left(2a^2c^2,a^2c^2(2ac - b^2)\right),\\
	[3]P=&\left(\frac{7}{9}a^2c^2 - \frac{16}{9}ab^2c + \frac{4}{9}b^4,-\frac{(2ac - b^2)(5a^2c^2 - 32ab^2c + 8b^4)}{27}\right).
	\end{split}
\end{equation*}
	We omit the expressions for the other six points because they will not be used directly. It is easy to verify that when $2ac-b^2\neq0$, the points $P$ and $\pm [i]P,~i=2,...,9,$ are different. By Mazur's theorem (see p. 58 of \cite{Silverman-Tate}) about the quantity of rational points and the rank of elliptic curve: If an elliptic curve $E$ defined over $\mathbb{Q}$ has more than 16 different rational points, then it has infinitely many rational points and its rank has at least one. Therefore, $\mathcal{E}_2$ has a positive rank, and thus there are infinitely many rational points on $\mathcal{E}_2$.
\end{prof}	

\begin{theorem}
	Let \(G=(x+y)z-2xy\) and \(f(x)=ax+b+cx^{-1} \in \mathbb{Q}[x,x^{-1}]\) with \(abc \neq 0\). When \(4ac - b^{2}\neq 0\), \(\langle G,F\rangle\) is birationally equivalent to the elliptic curve
	\[\mathcal{E}_3: Y^{2}=X^3-3a^2c^6X-a^2c^8(2ac-b^2).\]
When $4ac - b^{2}= 0$, \(\langle G,F\rangle\) is a curve with genus $0$ and its parameterization is given by
\[\begin{split}
y=\frac{8ct(9bt + 6ct + b)}{(17bt + 6ct + b)(5bt + 6ct + b)},~~z=-\frac{(17bt + 6ct + b)(9bt + 6ct + b)c}{16b^3t^2},
\end{split}\]
where $t$ is a rational number.
\end{theorem}

\begin{prof}
When \(G=(x+y)z-2xy\) and \(f(x)=ax+b+cx^{-1}\), the Diophantine system (\ref{1.2}) is equivalent to
\begin{equation}
	(y-z)^2\left(aby^2z^2+3acy^2z-2c^2y+c^2z\right)=0.                           \label{6.1}
\end{equation}
Since \(y=z\) is trivial, we only need to consider the rational points on the curve
\[\mathcal{C}_3: aby^2z^2+3acy^2z-2c^2y+c^2z=0.\]
By the map $\varphi_3$:
\begin{equation}
	\begin{split}
		X=&~\frac{c^2(aby^2z + 2acy^2 + c^2)}{y^2},\quad Y=~\frac{c^4(aby^3 + aby^2z + 3acy^2 + c^2)}{y^3},       \label{6.2}
	\end{split}
\end{equation}
with the inverse map $\varphi^{-1}_3$:
\begin{equation}
	\begin{split}
		y=&~\frac{c^2(Y+abc^4)}{(X+ac^3 )(X-2ac^3)},\quad  z=~\frac{2c^2(Y-abc^4)}{(X+ac^3)^2},                                                \label{6.3}
	\end{split}
\end{equation}
we can transform $\mathcal{C}_3$ into the elliptic curve
\begin{align*}
	\mathcal{E}_3:~Y^{2}=X^3 - 3a^2c^6X - a^2c^8(2ac - b^2).
\end{align*}
The discriminant of $\mathcal{E}_3$ is
\[\Delta_3(E)=27a^4c^{16}b^2(4ac - b^2).\]
Hence, if $4ac-b^2\neq0$, then $\Delta_3(E)\neq0$, so $\mathcal{E}_3$ is non-singular. Therefore, (\ref{6.2}) and (\ref{6.3}) give a bijection between \(\langle G,F\rangle\) and the elliptic curve $\mathcal{E}_3$, so \(\langle G,F\rangle\) is birationally equivalent to the elliptic curve $\mathcal{E}_3$.

When $4ac-b^2=0$, the Diophantine system (\ref{1.2}) is equivalent to
\begin{equation}
	b^3y^2z^2 + 3b^2cy^2z - 8c^3y + 4c^3z=0.                                    \label{6.4}
\end{equation}
The curve given by (\ref{6.4}) is a curve of genus $0$, and its parameterization is given by
\[\begin{split}
y=\frac{8ct(9bt + 6ct + b)}{(17bt + 6ct + b)(5bt + 6ct + b)},~~z=-\frac{(17bt + 6ct + b)(9bt + 6ct + b)c}{16b^3t^2},
\end{split}\]
where $t$ is a rational number. This completes the proof.
\end{prof}	

\begin{theorem}
	When \(b^{2}\neq kac~(k=\frac{27}{8},\frac{27}{4})\), the elliptic curve \(\mathcal{E}_3\) has a positive rank. When \(b^{2}=kac~(k=\frac{27}{8},\frac{27}{4})\), the rank of the elliptic curve \(\mathcal{E}_3\) is \(0\).
\end{theorem}

\begin{prof}
It is easy to check that the elliptic curve $\mathcal{E}_3$ contains the following rational point
	\begin{align*}
		P=(2ac^3, -abc^4).
	\end{align*}
By the group law, we get the following eight points
\begin{align*}
	&[2]P,\quad [3]P,\quad [4]P,\quad [5]P,\quad [6]P,\quad [7]P,\quad [8]P,\quad [9]P.
\end{align*}
The point $P_2$ is as follows
\begin{equation*}
	\begin{split}
		[2]P=&\left(\frac{ac^3(81ac - 16b^2)}{4b^2},\frac{c^4a(729a^2c^2 - 216ab^2c + 8b^4)}{8b^3}\right).
	\end{split}
\end{equation*}
We omit the expressions for the other seven points because they will not be used directly. It is easy to verify that when \(b^{2}\neq kac~(k=\frac{27}{8},\frac{27}{4})\), the points $P$ and $\pm [i]P,~i=2,...,9,$ are different. By Mazur's theorem (see p. 58 of \cite{Silverman-Tate}) about the quantity of rational points and the rank of elliptic curve: If an elliptic curve $E$ defined over $\mathbb{Q}$ has more than 16 different rational points, then it has infinitely many rational points and its rank has at least one. Therefore, $\mathcal{E}_3$ has a positive rank, and thus there are infinitely many rational points on $\mathcal{E}_3$.

When \(b^{2}=kac~(k=\frac{27}{8},\frac{27}{4})\), we have \(a=\frac{b^2}{ck}\), and then
\begin{align*}
	\mathcal{E}'_3:~Y^2=~X^3-\frac{3b^4c^4}{k^2}X+\frac{b^6c^6(k - 2)}{k^3}.
\end{align*}
By the transformation
\begin{equation}
	U=\frac{Y}{b^3c^3},\quad V=\frac{X}{b^2c^2},
\end{equation}
we get
\begin{align*}
	\mathcal{E}_{(k)}:~U^2=V^3-\frac{3}{k^2}V+\frac{k - 2}{k^3}.
\end{align*}

1) When $k=\frac{27}{8}$, we get
\begin{align*}
	\mathcal{E}_{(\frac{27}{8})}:~U^2=V^3-\frac{64}{243}V+\frac{704}{19683}.
\end{align*}
Using the package of Magma, the rank of $\mathcal{E}_{(\frac{27}{8})}$ is $0$.

2) When $k=\frac{27}{4}$, we get
\begin{align*}
    \mathcal{E}_{(\frac{27}{4})}:~U^2=V^3-\frac{16}{243}V+\frac{304}{19683}.
\end{align*}
The rank of $\mathcal{E}_{(\frac{27}{4})}$ is $0$.

Therefore, when \(b^{2}\neq kac~(k=\frac{27}{8},\frac{27}{4})\), the elliptic curve \(\mathcal{E}_3\) has a positive rank. When \(b^{2}=kac~(k=\frac{27}{8},\frac{27}{4})\), the rank of the elliptic curve \(\mathcal{E}_3\) is \(0\).
\end{prof}	

\begin{theorem}
	Let \(G=(x+y)z-2xy\) and \(f(x)=x(ax^2 + bx + c) \in \mathbb{Q}[x]\) with \(abc \neq 0\). When \(ac - b^{2}\neq 0, 3ac - b^{2}\neq 0\) and \(4ac - b^{2}\neq 0\), \(\langle G,F\rangle\) is birationally equivalent to the elliptic curve
	\[\mathcal{E}_4: Y^{2}=X^3 - 3a^6c^2X - a^6(2ac - b^2)(a^2c^2 - 4ab^2c + b^4).\]
When $ac - b^{2}=0$, \(\langle G,F\rangle\) is a curve with genus $0$ and its parameterization is given by
\[\begin{split}
y=-\frac{c(2t^2 + 6t + 3)}{3(t + 1)b},~~z=-\frac{c(2t^2 + 6t + 3)}{3b(t + 1)^2},
\end{split}\]
where $t$ is a rational number. When $3ac - b^{2}=0$, \(\langle G,F\rangle\) is the union of two curves, which are
\[z=-\frac{2c}{b},~~z=\frac{y(by + 2c)}{c}.\]
When $4ac - b^{2}=0$, \(\langle G,F\rangle\) is a curve with genus $0$ and its parameterization is given by
\[\begin{split}
y=-\frac{2ct(bt + 6ct + 3b)}{3(bt + 2ct + b)(2ct + b)},~~z=-\frac{4(2bt + 6ct + 3b)ct}{3(bt + 2ct + b)^2},
\end{split}\]
where $t$ is a rational number.
\end{theorem}

\begin{prof}
	When \(G=(x+y)z-2xy\) and \(f(x)=x(ax^2 + bx + c)\), the Diophantine system (\ref{1.2}) is equivalent to
	\begin{equation}
		(y-z)^2\left(3a^2y^2z + 2aby^2 + 2abyz - abz^2 - 2acy + acz + 2b^2y - b^2z\right)=0.                           \label{8.1}
	\end{equation}
Since \(y=z\) is trivial, we only need to consider the rational points on the curve
\[\mathcal{C}_4: 3a^2y^2z + 2aby^2 + 2abyz - abz^2 - 2acy + acz + 2b^2y - b^2z=0.\]
By the map $\varphi_4$:
\begin{equation}
	\begin{split}
		X=&~\frac{2a^3cy^2 - a^2b^2y^2 - a^2bcz + ab^3z + a^2c^2 - 2ab^2c + b^4}{y^2},\\
		Y=&~\frac{(ac - b^2)Y_1}{y^3},       \label{8.2}
	\end{split}
\end{equation}
where
\[\begin{split}
	Y_1=3a^3cy^2 - a^2b^2y^2 + a^2b^2yz - a^2bcy - a^2bcz+ ab^3y + ab^3z + a^2c^2 - 2ab^2c + b^4,
\end{split}\]
with the inverse map $\varphi^{-1}_4$:
\begin{equation}
	\begin{split}
		y=&~-\frac{(ac - b^2)(-2a^4bc + a^3b^3 + Xab - Y)}{(a^3c - 2a^2b^2 + X)(-2a^3c + a^2b^2 + X)},\\
		z=&~-\frac{2(ac - b^2)(-a^4bc - a^3b^3 + 2Xab - Y)}{(a^3c - 2a^2b^2 + X)^2},                                                \label{8.3}
	\end{split}
\end{equation}
we can transform $\mathcal{C}_4$ into the elliptic curve
\begin{align*}
	\mathcal{E}_4:~Y^{2}=X^3 - 3a^6c^2X - a^6(2ac - b^2)(a^2c^2 - 4ab^2c + b^4).
\end{align*}
The discriminant of $\mathcal{E}_4$ is
\[\Delta_4(E)=27a^{12}b^2(4ac - b^2)(3ac - b^2)^2(ac - b^2)^2.\]
Hence, if $ac - b^{2}\neq 0, 3ac - b^{2}\neq 0$ and $4ac-b^2\neq0$, then $\Delta_4(E)\neq0$, so $\mathcal{E}_4$ is non-singular. Therefore, (\ref{8.2}) and (\ref{8.3}) give a bijection between \(\langle G,F\rangle\) and the elliptic curve $\mathcal{E}_4$, so \(\langle G,F\rangle\) is birationally equivalent to the elliptic curve $\mathcal{E}_4$.

When $ac-b^2=0$, the Diophantine system (\ref{1.2}) is equivalent to
\begin{equation}
	3by^2z + 2cy^2 + 2cyz - cz^2=0.                                    \label{8.4}
\end{equation}
The curve given by (\ref{8.4}) is a curve of genus $0$, and its parameterization is given by
\[\begin{split}
y=-\frac{c(2t^2 + 6t + 3)}{3(t + 1)b},~~z=-\frac{c(2t^2 + 6t + 3)}{3b(t + 1)^2},
\end{split}\]
where $t$ is a rational number. When $3ac-b^2=0$, the Diophantine system (\ref{1.2}) is equivalent to
\begin{equation*}
	(bz + 2c)(by^2 + 2cy - cz)=0,
\end{equation*}
which leads to
\[z=-\frac{2c}{b},~~z=\frac{y(by + 2c)}{c}.\]
When $4ac-b^2=0$, the Diophantine system (\ref{1.2}) is equivalent to
\begin{equation}
	3b^2y^2z + 8bcy^2 + 8bcyz - 4bcz^2 + 24c^2y - 12c^2z=0.                                    \label{8.5}
\end{equation}
The curve given by (\ref{8.5}) is a curve of genus $0$, and its parameterization is given by
\[\begin{split}
y=-\frac{2ct(bt + 6ct + 3b)}{3(bt + 2ct + b)(2ct + b)},~~z=-\frac{4(2bt + 6ct + 3b)ct}{3(bt + 2ct + b)^2},
\end{split}\]
where $t$ is a rational number. This completes the proof.
\end{prof}	

\begin{theorem}
	When \(b^{2}\neq \frac32ac\), the elliptic curve \(\mathcal{E}_4\) has a positive rank. When \(b^{2}=\frac32ac\), the rank of the elliptic curve \(\mathcal{E}_4\) is \(0\).
\end{theorem}

\begin{prof}
It is easy to check that the elliptic curve $\mathcal{E}_4$ contains two rational points
\begin{align*}
	P_0=(2a^3c - a^2b^2,0),~~P_1=(-a^3c, a^3b(3ac - b^2)).
\end{align*}
By the group law, we get the following seven points
\begin{align*}
	&P_2=[2]P_1,\quad P_3=[3]P_1,\quad P_4=[4]P_1,\\
	P_5=P_1+&P_0,\quad P_6=P_2+P_0,\quad P_7=P_3+P_0,\quad P_8=P_4+P_0.
\end{align*}
The point $P_2$ is as follows
\begin{align*}
	P_2=(2a^3c,-3cba^4 + a^3b^3).
\end{align*}
We omit the expressions for the other six points because they will not be used directly. It is easy to verify that when \(b^{2}\neq \frac32ac\), the points $P$ and $\pm [i]P,~i=2,...,9,$ are different. By Mazur's theorem (see p. 58 of \cite{Silverman-Tate}) about the quantity of rational points and the rank of elliptic curve: If an elliptic curve $E$ defined over $\mathbb{Q}$ has more than 16 different rational points, then it has infinitely many rational points and its rank has at least one. Therefore, $\mathcal{E}_4$ has a positive rank, and thus there are infinitely many rational points on $\mathcal{E}_4$.

When \(b^{2}=\frac32ac\), we have \(a=\frac{2b^2}{3c}\), and then
\begin{align*}
	\mathcal{E}'_4:~Y^2=~X^3-\frac{64b^{12}}{243c^4}X+\frac{704b^{18}}{19683c^6}.
\end{align*}
By the transformation
\begin{equation}
	U=\frac{c^3Y}{b^9},\quad V=\frac{c^2X}{b^6},
\end{equation}
we get
\begin{align*}
\mathcal{E}''_4:~U^2=V^3-\frac{64}{243}V+\frac{704}{19683}.
\end{align*}
Using the package of Magma, the rank of $\mathcal{E}''_4$ is $0$.
\end{prof}

\end{document}